\documentclass[11pt]{article}

\catcode`\@=11
\def\marginnote#1{}

\newcount\hour
\newcount\minute
\newtoks\amorpm
\hour=\time\divide\hour by60 \minute=\time{\multiply\hour by60
\global\advance\minute by-\hour}
\edef\standardtime{{\ifnum\hour<12
\global\amorpm={am}%
        \else\global\amorpm={pm}\advance\hour by-12 \fi
        \ifnum\hour=0 \hour=12 \fi
        \number\hour:\ifnum\minute<10
        0\fi\number\minute\the\amorpm}}
\edef\militarytime{\number\hour:\ifnum\minute<10
0\fi\number\minute}

\def\draftlabel#1{{\@bsphack\if@filesw {\let\thepage\relax
   \xdef\@gtempa{\write\@auxout{\string
      \newlabel{#1}{{\@currentlabel}{\thepage}}}}}\@gtempa
   \if@nobreak \ifvmode\nobreak\fi\fi\fi\@esphack}
        \gdef\@eqnlabel{#1}}
\def\@eqnlabel{}
\def\@vacuum{}
\def\draftmarginnote#1{\marginpar{\raggedright\scriptsize\tt#1}}

\def\draft{\oddsidemargin -.5truein
        \def\@oddfoot{\sl preliminary draft \hfil
        \rm\thepage\hfil\sl\today\quad\militarytime}
        \let\@evenfoot\@oddfoot \overfullrule 3pt
        \let\label=\draftlabel
        \let\marginnote=\draftmarginnote
   \def\@eqnnum{(\theequation)
   \rlap{\kern\marginparsep\tt\@eqnlabel}%
\global\let\@eqnlabel\@vacuum} } 
\relax

\draft
\input{amssym.def}
\input{amssym}
\usepackage{eufrak}
\usepackage{color}

\newtheorem{theorem}{Theorem}

\newtheorem{proposition}[theorem]{Proposition}
\newtheorem{remark}[theorem]{Remark}

\def\pa{\partial}
\def\h{h}
\def\rh{r_\hbar}
\def\De{\Delta}
\def\D{\cal D}
\def\Cas{{\rm Cas}}

\def\tpa{\tilde{\pa}_t}

\def\AA{\cal A}
\def\E{{\bf  E}}
\def\H{{\bf H}}

\def\di{\rm div}
\def\Ah{{\frak{A}}}
\def\AAA{{{\bf{A}}}}
\def\hh{\hbar}

\def\nn{\overrightarrow{n}}

\def\parh{\pa_{\rh}}
\def\tpa{\tilde{\pa}}
\def\UU{U(gl(m)_h)}
\def\Mat{\rm Mat}

\def\A{\cal{A}}
\def\de{\delta}

\def\gg{\mbox{$\frak g$}}
\def\ot{\otimes}

\def\K{{\Bbb K}}

\def\C{{\Bbb C}}
\def\R{{\Bbb R}}
\def\AA{{\cal A}}
\def\UU{{U(u(2)_h)}}
\def\Sym{{\rm Sym\, }}

\def\End{{\rm End}}

\def\lhq{\ifmmode {\cal L}(q,\hbar)\else ${\cal L}(q,\hbar)$\fi}

\def\lqh{\ifmmode {\cal L}(q,\hbar)\else ${\cal L}(q,\hbar)$\fi}

\def\W{{\cal W}}
\def\gh{\gg_{h}}

\def\al{{\alpha}}
\def\hatt{\hat{\Theta}}

\def\be{\begin{equation}}
\def\ee{\end{equation}}

\begin{document}

\makeatletter
\renewcommand{\theequation}{{\thesection}.{\arabic{equation}}}
\@addtoreset{equation}{section} \makeatother

\title{Quantum geometry and  quantization on $U(u(2))$ background.\\
Noncommutative Dirac monopole}

\author{\rule{0pt}{7mm} Dimitri
Gurevich\thanks{gurevich@ihes.fr}\\
{\small\it LAMAV, Universit\'e de Valenciennes,
59313 Valenciennes, France}\\
\rule{0pt}{7mm} Pavel Saponov\thanks{Pavel.Saponov@ihep.ru}\\
{\small\it
National Research University Higher School of Economics,}\\
{\small \it International Laboratory of Representation Theory and  Mathematical Physics}\\
{\small\it 20 Myasnitskaya Ulitsa, Moscow 101000, Russia}\\
{\small\it and}\\
{\small\it
State Research Center of the Russian Federation  ``Institute for High Energy Physic''}\\
{\small\it of National Research Centre ``Kurchatov Institute''}\\
{\small\it Protvino, Russia}}

\maketitle

\begin{abstract}
In our previous publications we introduced differential calculus on the enveloping algebras  $U(gl(m))$ similar to the usual calculus on the commutative algebra 
$\Sym\!(gl(m))$. The main ingredient of our calculus are quantum partial derivatives which turn into the usual partial derivatives in the classical limit. In the particular 
case $m=2$ we prolonged this calculus on a central extension $\AA$ of the algebra $U(gl(2))$. In the present paper we consider the problem of a further 
extension of the quantum partial derivatives on the skew-field of the algebra $\AA$ and define the corresponding de Rham complex. As an application of the 
differential calculus we suggest a method of transferring dynamical models defined on ${\Sym}\!(u(2))$ to the algebra $U(u(2))$ (we call this procedure the 
quantization with noncommutative configuration space). In this sense we quantize the Dirac monopole and find a solution of this model.
\end{abstract}

{\bf AMS Mathematics Subject Classification, 2010:} 81S99, 83C65


\section{Introduction}

By Quantum Geometry  we mean a sort of Noncommutative (NC) Geometry dealing with objects and operators which are deformations of the classical ones.
The main ingredients of this type Geometry are analogs of partial derivatives defined on the algebras\footnote{Given a Lie algebra $\gg$ with the Lie bracket 
$[\,,\,]$, the notation $\gh$ stands for the Lie algebra whose bracket  is $h\,[\,,\,]$, where $h$ is a quantization parameter. Note that the parameter $h$ is not 
formal and can be specialised.} $U(gl(m)_h)$ or their extentions. These {\em quantum partial derivatives} have the mentioned deformation property: for $\h=0$ 
they turn into the usual partial derivatives on the algebra $\Sym\!(gl(m))$. Besides, similarly to the usual partial derivatives, the quantum\footnote{Below, we 
often omit the precision {\em quantum}.} ones commute with each other.

This Quantum Geometry enables us to quantize dynamical models in the following sense. Consider a differential operator, corresponding to such a model. If 
the coefficients of the operator belong to the algebra $\Sym\!(gl(m))$, we replace them by their noncommutative counterparts from 
the algebra $U(gl(m)_h)$. Then, on replacing the usual partial derivatives by their quantum counterparts, we get a quantum differential operator well defined on the
algebra $U(gl(m)_h)$. We call this procedure the quantization with {\em noncommutative configuration space} since the final operators acts on a NC algebra and 
not on a Hilbert space.

However, it is not clear how to quantize an operator with more general coefficients, e.g. rational functions. In order to quantize such an operator, we have to extend
the quantization of the algebra $\Sym\!(gl(m))$ to the field of rational functions, as well as to extend the quantum partial derivatives to the skew-field (division ring) 
of the algebra $U(gl(m)_h)$. A similar question arises while we add some elements of an algebraic central extension to the initial algebra. In the present
paper we give an answer to this question in a particular case $m=2$ (in fact, we deal with the compact form $\UU$ of the algebra $U(gl(2)_h)$).
The ground field $\K$ is assumed to be $\C$ (except the last section where $\K=\R$).

Now, we briefly recall a way of introducing the quantum partial derivatives. The way, which is the most similar to the classical one, is based on the coproduct defined 
on the commutative algebra $\D$, which is generated by the quantum partial derivatives. Let $\{n_k^l\}$ be the standard basis in the algebra $U(gl(m)_h)$. Define an
action of a quantum partial derivative $\pa_j^i=\pa_{n_i^j}$ on the generators $n_k^s$ by the classical formula $\pa_i^j(n_k^s)=\de_i^s\de_k^j$ and extend this
action on higher order monomials in the generators by means of the following coproduct:
\be
\De(\pa_i^j)= \pa_i^j\ot 1+1\ot  \pa_i^j +\frac{\h}{2} \sum_k \pa^j_k \ot \pa_i^k. 
\label{co} 
\ee
Note that this coproduct is coassociative  but not cocommutative.

Emphasize  that for $\h=0$ we get the standard coproduct for the classical partial derivatives on the algebra $\Sym\!(gl(m))$ (the usual Leibniz rule). The coproduct
(\ref{co}) enables us to compute the action of a quantum partial derivative $\pa_i^j$ on a product $ab$, provided we know the action of {\it all} partial derivatives on 
$a\in U(gl(m)_h)$ and $b\in U(gl(m)_h)$ separately. Whereas, in the classical case it suffices to know the action of {\it the same} derivative $\pa_i^j$ on the factors.

This coproduct also enables us to define the so-called {\em permutation relations} in the following way:
\be
\pa\ot a \mapsto \pa_{(1)}(a)\ot \pa_{(2)},\quad \forall \pa\in {\D},\,\, \forall a\in U(gl(m)_h),
\label{first} 
\ee
where $\De(\pa)=\pa_{(1)}\ot \pa_{(2)}$ in the Sweedler's notation.

Thus, formula (\ref{first}) gives rise to a map
\be
\sigma: {\D}\ot U(gl(m)_h)\to U(gl(m)_h)\ot \D, 
\label{map} 
\ee
enabling us to introduce an associative product in $U(gl(m)_h)\ot \D$ by the formula
$$
(a_1\ot \pa)(a_2\ot \pa')=(\cdot \ot \cdot)(I\ot \sigma\ot I)(a_1\ot \pa\ot a_2\ot \pa'),\qquad a_1,a_2\in U(gl(m)_h),\quad \pa,\pa'\in \D.
$$
Here $\cdot\ot \cdot$ stands for products in the factors.

The resulting algebra $U(gl(m)_h)\ot \D$ is called {\em the quantum Weyl-Heisenberg (WH) algebra}\footnote{Mathematicians call this algebra Weyl one, physicists call
it Heisenberg one. We proceed here in the nature of compromise.} and is denoted $\W(U(gl(m)_h))$. Thus, the main role in the construction of the WH algebra is played
by the coproduct (\ref{co}) which is a deformation of the classical one. However, there exist algebras, related to {\em braidings} and close in a sense to the enveloping
algebras $U(gl(m)_h)$, in which the quantum partial derivatives and analogs of the WH algebras can also be introduced, but no coproduct similar to (\ref{co}) is known.  
We do not consider such  algebras here and refer the reader to the paper \cite{GS1, GS3}, where they are studied and the corresponding differential calculus is constructed. Thus, the permutation relations constitute the second (and in principle, more general) way for introducing quantum partial derivatives.

Once the partial derivatives are defined, we can construct an analog of the de Rham complex with the usual property $d^2=0$ of the de Rham operator. This complex 
differs from the Hochschild or cyclic ones, usually employed in other approaches to NC Geometry. In contrast  with the latter complexes, the terms of ours are deformations of
their classical counterparts.

In the current paper we deal with the algebra $\UU$ (which can be viewed as a deformed space-time algebra). The main objective of this paper is an extension of our differential
calculus on $\UU$ up to a bigger algebra, which would include some central algebraic extension and the division ring of the algebra $\UU$. The main problem consists in finding
a proper extension of the partial derivatives on the mentioned algebra. In order to solve this problem, we introduce a special matrix $\hatt(a)$ composed of the partial derivatives
and consider a map $\UU\ni a\to \hatt(a)$, where the notation $\hatt(a)$ means that any entry of the matrix $\hatt(a)$ is applied to a given element $a$. This map is multiplicative. 
We treat this property as one of the forms of the Leibniz rule. In order to define an action of the partial derivatives on the element $a^{-1}$, $a\not=0$ we have to invert the matrix
$\hatt(a)$. After having reviewed different forms of the Leibniz rule for the partial derivatives on the algebra $\UU$ (section 2), we consider the problem of explicit inverting the matrices $\hatt(a)$ in function of a given element $a$. Finally, in section 3 we extend the partial derivatives and other elements of our differential calculus on a family containing
the algebra $\UU$. In section 4 we explain how this differential calculus can be used in a quantization of dynamical models in a new sense, which consists in transferring these 
models to noncommutative configuration space. In section 5 we apply this method for quantization (in our sense) of the Dirac monopole. As a result we arrive to a NC version of
this model.

\medskip

\noindent
{\bf Acknowledgement.} 
This research was supported jointly by the National Research
University ---``Higher School of Economics'' Academic Fund Program (grant
No.15-09-0279) and by the grant of RFBR No.14-01-00474-a.

\section{Quantum partial derivatives: different forms of the Leibniz rule}

In this section we consider different ways of introducing the quantum partial derivatives and corresponding forms of the Leibniz rule. First, we recall the definition of the 
Weyl-Heisenberg algebra $\W(U(u(2)_h))$. Similarly to the classical case, it is generated by two subalgebras. One of them is the enveloping algebra $U(u(2)_h)$. It is 
generated by the unit $1_{U(u(2)_h)}$ and elements $\{t, x,y,z\}$ subject to the relations
$$
[x, \, y]= \h\, z,\quad [y, \, z]= \h\, x,\quad[z, \, x]=\h\, y,\quad [t, \, x]=[t, \, y]=[t, \, z]=0.
$$
Note that the element $t$ and the quadratic Casimir element
$$
{{\Cas}} = x^2+y^2+z^2
$$
are central. Moreover, the center $Z(U(u(2)_h))$ of the algebra $U(u(2)_h)$ is generated by these two elements.

The second subalgebra of the algebra $\W(U(u(2)_h))$ (denoted $\D$) is generated by the unit $1_{\D}$ and by the commutative elements  $\partial_t$, $\partial_x$, 
$\partial_y$, $\partial_z$. Also, we impose the following {\em permutation relations} among the generators 
$\partial_t$, $\partial_x$, $\partial_y$, $\partial_z\in\D$ and $t, x,y,z\in\UU$:
\be
\begin{array}{l@{\quad}l@{\quad}l@{\quad}l}
\tilde\pa_t\,t - t\,\tilde\pa_t = {h\over 2}\,\tilde\pa_t & \tilde\pa_t\, x - x\,\tilde\pa_t
=-{h\over 2}\,\pa_x &
\tilde\pa_t\, y - y\, \tilde\pa_t=-{h\over 2}\,\pa_y &\tilde\pa_t\, z - z\,\tilde\pa_t=- {h\over 2}\,\pa_z\\
\rule{0pt}{7mm}
\pa_x\, t - t\,\pa_x = {h\over 2}\,\pa_x &\pa_x \,x -  x\,\pa_x = {h\over 2}\,\tilde\pa_t &
\pa_x \, y-  y\,\pa_x = {h\over 2}\,\pa_z & \pa_x \,z - z\, \pa_x  = - {h\over 2}\,\pa_y \\
\rule{0pt}{7mm}
\pa_y \,t - t \, \pa_y = {h\over 2}\,\pa_y & \pa_y \,x -  x\,  \pa_y = -{h\over 2}\,\pa_z &
\pa_y \,y - y \,  \pa_y = {h\over 2}\,\tilde\pa_t & \pa_y \,z - z \,  \pa_ y= {h\over 2}\,\pa_x\\
\rule{0pt}{7mm}
\pa_z \,t - t \,\pa_z = {h\over 2}\,\pa_z & \pa_z \,x - x \,\pa_z = {h\over 2}\,\pa_y&
\pa_z \,y -  y\,\pa_z = -{h\over 2}\,\pa_x & \pa_z \,z - z \,\pa_z = {h\over 2}\,\tilde\pa_t.
\end{array}
\label{perm}
\ee

Hereafter, the notation $\tpa_t$ stands for the {\em shifted derivative} in $t$:
$$\tilde \partial_t = \partial_t +\frac{2}{h}1_{\D}.$$

Also, we assume that the unit $1_{U(u(2)_h)}$ (resp., $1_{\D}$) permutes with elements
of $\D$ (resp., $U(u(2)_h)$) in the standard  way, i.e. by means of the usual flip.
Extending the above permutation relations to the whole algebras $U(u(2)_h)$ and $\D$, we get the map $\sigma$   (\ref{map}).
Note, that the elements $\pa_u,\,\,u\in\{x,y,z\} $ in the r.h.s. of (\ref{perm}) must be treated as $1_{U(u(2)_h)}\ot \pa_u$, and the 
term $\tpa_t=\pa_t+\frac{2}{h}$ as $1_{U(u(2)_h)}\ot \tpa_t $.
Hereafter,  we omit the symbol $\ot$ when it does not lead to a misunderstanding.

In order to convert the quantum partial derivatives into operators we need the counit $\varepsilon:\D\to \K$ which is an algebra homomorphism defined on 
the generators by
\be
\varepsilon(1_{\D})=1, \quad \varepsilon(\pa_u)=0,\quad \forall u\in \{t,x,y,z\}. 
\label{coun} 
\ee
Taking into account the definition of $\tilde\partial_t$, we get $\varepsilon(\tpa_t)={{2}\over {h}}$.
Then, in order to find an action of a partial derivative $\pa_u$ on an element $a\in  U(u(2)_h)$ we permute the elements $\pa_u$ and $a$ in the product
$\pa_u\ot a$ with the use of the relations (\ref{perm}) and apply the above counit map to the right factor, belonging to $\D$. The resulting element is denoted $\pa_u(a)$.

This definition immediately entails that $\pa_x(x)=1$, $\pa_x(y)= \pa_x(z)=0$ and so on, i.e. the partial derivatives act on the generators of the algebra $U(u(2)_h)$
is  the same way as in the algebra $\Sym\!(u(2))$.  By contrary, the action of the partial derivatives on higher degree elements of the algebra $U(u(2)_h)$ differ from 
the classical one.

Let us consider an example. Using the permutation relations, we get
$$
\pa_x yz=(y\pa_x+\frac{h}{2}\pa_z)z=y(z\pa_x-\frac{h}{2}\pa_y)+\frac{h}{2}(z\pa_z+\frac{h}{2}\tpa_t)=
yz\pa_x-\frac{h}{2}y\pa_y+\frac{h}{2}z\pa_z+\frac{h^2}{4}\tpa_t.
$$
Now, by applying the counit to factors belonging to the algebra $\D$, we get $\pa_x (yz)=\frac{h}{2}$.
In the same way we get $\pa_x (zy)=-\frac{h}{2}$. This result is compatible with the relation  $yz-zy=hx$.

Thus, the permutation relations (\ref{perm}) together with the counit $\varepsilon $ can be considered as an analog of the classical Leibniz rule. Another form the Leibniz 
rule for the partial derivatives on the algebra $U(u(2)_h)$ can be presented as follows.
For $h=1$ we identify the generators $t$, $x$, $y$, $z$ with elements of the algebra $\End(V)$, where   $V$ is the fundamental $u(2)$-module. Namely, we have
\be
 t=\frac{1}{2}\left(\begin{array}{cc}
1&0\\
0&1
\end{array}\right),\quad x=\frac{i}{2}\left(\begin{array}{cc}
0&1\\
1&0
\end{array}\right),\quad y=\frac{1}{2}\left(\begin{array}{cc}
0&-1\\
1&0
\end{array}\right),\quad z=\frac{i}{2}\left(\begin{array}{cc}
1&0\\
0&-1
\end{array}\right).
\label{pau}
\ee

Let $\circ: \End(V)^{\ot 2}\to \End(V)$ be the usual product, namely, the composition of endomorphisms. Then we have the following multiplication table for the 
elements $t,x,y,z$:
$$
t\circ u=u\circ t=\frac{u}{2},\quad u\in\{t, x,y,z\},\quad x\circ x=-\frac{t}{2},\quad x\circ y=-y\circ x =\frac{z}{2}, \,\,\circlearrowleft
$$
where $\circlearrowleft$ stand for the cyclic permutations $x\to y \to z$. Note, that the $u(2)$ Lie bracket is related to this product as follows
$[u,v]=u\circ v-v\circ u$, $\forall \,  u,v\in \End(V)$.

Then the action  $\pa_u(ab)$, $u$, $a$, $b \in \{t, x,y,z\}$  is defined by
$$
\pa_u(a b)= \pa_u(a)b+a\pa_u(b)+h\pa_u(a\circ b).
$$
This rule together with its extension to the higher degree polynomials was called in \cite{GS2} the $h$-Leibniz rule.

One more form of the Leibniz rule is based on a coproduct
$$
\De:U(u(2)_h)\to U(u(2)_h)\ot U(u(2)_h),
$$
 which in fact is a particular case of that (\ref{co}) realized in
 the basis $\{\pa_t, \pa_x,\pa_y,\pa_z\}$. In the explicit form it reads:
 \be
 \begin{array}{l}
\De(\pa_t)=\pa_t\ot 1 +1\ot \pa_t+\frac{h}{2}(\pa_t\ot \pa_t-\pa_x\ot \pa_x-\pa_y\ot \pa_y-\pa_z\ot \pa_z),\\
\rule{0pt}{5mm}
\De(\pa_x)=\pa_x\ot 1 +1\ot \pa_x+\frac{h}{2}(\pa_t\ot \pa_x+\pa_x\ot \pa_t+\pa_y\ot \pa_z-\pa_z\ot \pa_y),\\
\rule{0pt}{5mm}
\De(\pa_y)=\pa_y\ot 1 +1\ot \pa_y+\frac{h}{2}(\pa_t\ot \pa_y+\pa_y\ot \pa_t+\pa_z\ot \pa_x-\pa_x\ot \pa_z),\\
\rule{0pt}{5mm}
\De(\pa_z)=\pa_z\ot 1 +1\ot \pa_z+\frac{h}{2}(\pa_t\ot \pa_z+\pa_z\ot \pa_t+\pa_x\ot \pa_y-\pa_y\ot \pa_x).
\end{array} \label{arr1}
\ee

Being rewritten in terms of the shifted derivative $\tpa_t$, the coproduct (\ref{arr1}) takes the form
\be
\begin{array}{l}
\De(\tpa_t)=\frac{h}{2}(\tpa_t\ot \tpa_t-\pa_x\ot \pa_x-\pa_y\ot \pa_y-\pa_z\ot \pa_z),\\
\rule{0pt}{5mm}
\rule{0pt}{5mm}
\De(\pa_x)=\frac{h}{2}(\tpa_t\ot \pa_x+\pa_x\ot \tpa_t+\pa_y\ot \pa_z-\pa_z\ot \pa_y),\\
\rule{0pt}{5mm}
\De(\pa_y)=\frac{h}{2}(\tpa_t\ot \pa_y+\pa_y\ot \tpa_t+\pa_z\ot \pa_x-\pa_x\ot \pa_z),\\
\rule{0pt}{5mm}
\De(\pa_z)=\frac{h}{2}(\tpa_t\ot \pa_z+\pa_z\ot \tpa_t+\pa_x\ot \pa_y-\pa_y\ot \pa_x).
\end{array} \label{arr2} \ee

So, the shifted derivative $\tpa_t$ allows one to present this coproduct  in the multiplicative  form (\ref{arr2}) whereas its initial form (\ref{arr1}) is
additive-multiplicative.

Now, consider the permutation relations between the column
\be
(\tpa_t,\pa_x, \pa_y,\pa_z)^T \label{col} 
\ee
($T$ stands for the transposition) and any element $a\in U(u(2)_h)$. We have
\be
\left(\begin{array}{l}
\tpa_t\\
\pa_x\\
\pa_y\\
\pa_z\end{array}\right)\,a=\frac{h}{2}{\Theta}(a)
\left(\begin{array}{l}
\tpa_t\\
\pa_x\\
\pa_y\\
\pa_z\end{array}\right),\,\,{\rm where}\,\,
{\Theta}=\left(\begin{array}{llll}
\tpa_t&-\pa_x&-\pa_y&-\pa_z\\
\pa_x&\,\,\,\,\,\tpa_t&-\pa_z&\,\,\,\,\,\pa_y\\
\pa_y&\,\,\,\,\,\pa_z&\,\,\,\,\,\tpa_t&-\pa_x\\
\pa_z&-\pa_y&\,\,\,\,\,\pa_x&\,\,\,\,\,\tpa_t \end{array} \right). \label{seven}
\ee
This is a direct consequence of (\ref{arr2}) and (\ref{first}). Here, the notation ${\Theta}(a)$ stands for the matrix whose entries result from applying 
the corresponding entries of the matrix ${\Theta}$ to  $a$.

Also, note that the coproduct $\De$, being applied to the matrix ${\Theta}$ (i.e. to each entry),
leads to the formula
$$
\De ({\Theta})=\frac{h}{2}\, {\Theta} \stackrel{.}{\otimes} {\Theta},
$$
where the notation $A\stackrel{.}{\otimes} B$ stands for the matrix with entries $(A\stackrel{.}{\otimes} B)_i^j=\sum_k\, A_i^k\ot B_k^j $.

This entails the following relation
$$ 
{\Theta}(ab)=\frac{h}{2} {\Theta}(a){\Theta}(b)=i\hh {\Theta}(a){\Theta}(b),\quad \forall\, a,b \in \UU. 
$$
Hereafter, we also use another parameter $\hh$ which differs from $h$ by a factor: $h=2i\hh$.

Thus, up to a  factor the map $\UU\ni a \to \Theta(a)$ is a morphism of the algebra $\UU$ into the algebra ${\Mat}(\UU)$ of $4\times 4$ matrices
with coefficients belonging to $\UU$. More precisely, the map $\hatt=i\hh \Theta$ is a morphism:
\be 
{\hatt}(ab)={\hatt}(a){\hatt}(b),\quad \forall \,a,b \in \UU. 
\label{The}
\ee
In particular, $\hatt(1_\UU)=I$. Hereafter, by abusing the notation, we treat $\Theta$ and $\hatt$  as matrices and as the corresponding maps 
from $\UU$ to ${\Mat}(\UU)$.

Let us exhibit the images of the elements $x,y,z$, and $\Cas$ under the map $\hatt$:
\be
{\hatt}(x)=
\left(\!\!
\begin{array}{cccc}
 x &    -i\hh\,&0&0\\
  i\hh&x&0&0\\
0&0&  x & -i\hh\\
0&0&  i\hh & x
\end{array}
\!\!\right),
{\hatt}(y)=
\left(\!\!
\begin{array}{cccc}
 y &    0\,&-i\hh&0\\
0 &y&0&i\hh\\
i\hh&0&  y & 0\\
0&-i\hh&  0 & y
\end{array}
\!\!\right),
{\hatt}(z)=
\left(\!\!
\begin{array}{cccc}
 z & 0&0&-i\hh\\
 0 &z&-i\hh&0\\
0&i\hh&  z& 0\\
i\hh&0& 0 & z
\end{array}
\!\!\right),
\ee
\be
 {\hatt}(\Cas)=
\left(\!\!
\begin{array}{cccc}
\Cas+3\hh^2 &   \displaystyle -2i\hh x\,& \displaystyle -2i\hh y &\displaystyle -2i\hh z\\
\displaystyle  2i\hh x&\Cas+3\hh^2&\displaystyle -2i\hh z &\displaystyle 2i\hh y \\
\displaystyle 2i\hh y &\displaystyle 2i\hh z &\Cas+3\hh^2 &\displaystyle -2i\hh x\\
\displaystyle 2i\hh z &\displaystyle -2i\hh y & \displaystyle 2i\hh x&\Cas+3\hh^2
\end{array}
\!\!\right). \ee

Also, note that  $\hatt(t^p)=(t+i\hh)^p I,\,\, \forall p\in \Bbb{Z }$.

We treat the multiplicativity of the map $\hatt$ as another form of the Leibniz rule.
In the next section we define an extension of the algebra $U(u(2)_h)$ and
prolong the partial derivatives and consequently the map $\hatt$ on the extended algebra so that this multiplicativity remains valid.
In this sense we will speak about an extension of the WH algebra $\W(\UU)$.

\section{Extending  quantum partial derivatives}

The mentioned extension of the algebra $\UU$ will be constructed in two steps.
First, construct  a central algebraic extension of this algebra defined as follows.
Consider the matrix composed of the generators of the algebra $U(u(2)_h)$
\be
N=\left(\begin{array}{cc}
t-i\,z& -i\, x-y\\
-i\, x+y& t+i\, z
\end{array}\right).
\label{matrN}
\ee

\begin{remark} \rm This matrix is useful for   realizing the defining relations of the algebra $U(u(2)_h))$ in a concise  form:
$$
PN_1PN_1-N_1PN_1P=\h(PN_1-N_1P),
$$
where $P$ is the matrix of the usual flip $u\ot v \to  v\ot u$ and $N_1= N\otimes I$.
\end{remark}

The matrix $N$, called {\em the generating} matrix of the algebra $U(u(2)_h))$,
 is subject to the following NC analog of the Cayley-Hamilton (CH) identity:
\be
\chi(N) = N^2-(2\,t+h)\, N+\,(t^2+x^2+y^2+z^2+h\, t)\,I= 0.
\label{CH}
\ee

Note  that the coefficients  of the polynomial $\chi(N)$, are scalar
and belong to  the center $Z(U(u(2)_h))$ of the algebra $U(u(2)_h)$.

\begin{remark} 
\rm Note that a CH identity with similar properties is valid for  generating  matrices of the algebras $U(gl(m)_h)$ and their super- and braided analogs
(see \cite{GPS1}). In all these cases the generating matrix  is composed of generators of the corresponding algebra in a special way. In \cite{GPS1} a 
version of the CH identity was also suggested for generating matrices of other quantum matrix algebras, in particular the RTT one. However, for them the
powers of such a generating matrix are not treated in the usual sense and the coefficients of the corresponding polynomial $\chi$ are not central.
 \end{remark}

The roots $\mu_1$ and $\mu_2$ of the  polynomial $\chi(\mu)$ are called the {\em eigenvalues} of the matrix $N$. They are elements of an algebraic 
extension of the center  $Z(U(u(2)_h))$.

By expressing the generators $t$ and ${\Cas}$ of the center $Z(U(u(2)_h))$ via the eigenvalues  $\mu_1$ and $\mu_2$ we get
\be
t=\frac{\mu_1+\mu_2-\h}{2},\qquad {{\Cas}}=\frac{\h^2-\mu^2}{4}.
\label{cas}
\ee
where $\mu=\mu_1-\mu_2$. Below we do not need the eigenvalues and only include the quantity $\mu$ in the extended algebra, where  $\mu$ is assumed to be central.

More precisely, we deal with the following quantity
\be
\rh=\frac{\mu}{2i}=\sqrt{\Cas+\hh^2} \label{sign} 
\ee
which is called {\em the quantum radius} (recall that $h=2i\hh$). Let $\K(t,\rh) $ stand for the field of all rational functions in $t$ and $\rh$ (if necessary, $\K(t,\rh)$
can be enlarged up to the field of meromorphic functions or formal series).  Now, we put
\be
\AA=\left(U(su(2)_h)\ot \K(t,\rh)\right)/\langle x^2+y^2+z^2-\rh^2+\hh^2 \rangle. 
\label{quot} 
\ee
Hereafter, $\langle J \rangle$ stands for the ideal generated by a subset $J$.

The quantum radius is central in the  algebra $\AA$, since it is so for the quantity $\mu$.
In order to fix the sign of the square root in (\ref{sign}), we assume  the quantum
radius to be positive provided $\hh$ is real and $x,y,z$ are  represented by Hermitian operators.

In \cite{GS2} we extended the quantum partial derivatives to the algebra $\AA$. This enables us to compute the image $\hatt(a)$
of any element  $a\in \AA$. For instant, let us exhibit the matrix $\hatt(\rh)$:
\be {\hatt}(\rh)=
\frac{1}{\rh}\left(\!\!
\begin{array}{cccc}
\rh^2+\hh^2 &   \displaystyle -i\hh x\,& \displaystyle -i\hh y &\displaystyle -i\hh z\\
\displaystyle  i\hh x&\rh^2+\hh^2&\displaystyle -i\hh z &\displaystyle i\hh y \\
\displaystyle i\hh y &\displaystyle i\hh z & \rh^2+\hh^2 &\displaystyle -i\hh x\\
\displaystyle i\hh z &\displaystyle -i\hh y & \displaystyle i\hh x& \rh^2+\hh^2
\end{array}
\!\!\right).
 \label{mat1} \ee

We extend the action of the map $\hatt$ onto the algebra $\K(t,\rh)$ in a natural way by setting
$$
\hatt(f(t,\rh))=f(\hatt(t), \hatt(\rh)),\quad f\in \K(t,\rh).
$$
\begin{proposition} 
The map $\hatt: \AA\to {\Mat}(\AA)$ is well defined, i.e. the matrix $\hatt(\rh)$ commutes with any matrix $\hatt(a),\, a\in \UU$
and the map $\hatt$ is compatible with the relation $\Cas=\rh^2-\hh^2$, i.e.
\be
\hatt(\Cas)=\hatt(\rh)^2-\hh^2 I. 
\label{ree} 
\ee
\end{proposition}

{\bf Proof.}
The both claims can be proved by straightforward computations. Note, that it suffices to verify the first claim for $a\in \{x,y,z\}$ only.

\begin{remark} \rm In fact, the relation (\ref{ree}) can be used in order to find the matrix ${\hatt}(\rh)$. In \cite{GS1, GS2} we employed a similar method for finding
the result of applying the partial derivatives to $\rh^p$ but we dealt with other matrices and their spectral decompositions.
\end{remark}

Now, pass to the second step of extending the algebra $\UU$ and consider the skew-field $\Ah=\AA[\AA^{-1}]$. This field consists of left fractions 
$a^{-1} b$, $a,b \in\AA$, $a\not=0$ which in virtue of the Ore property  can be presented as right fractions $c\,d^{-1}$. We would like to extend the map $\hatt$ to 
the skew-field  $\Ah$, but in fact, we will extend the map $\hatt$ only on its subset.

First, recall the method of solving the same problem in a commutative unital algebra $A$.
Let  $V:A\to A$ be a vector field,  i.e. an operator  subject to the usual Leibniz rule.
Then applying $V$ to the element $a\,a^{-1}=1$, we get
$$
V(a)\,a^{-1}+a\,V(a^{-1})=0\quad \Rightarrow\quad V(a^{-1})=-V(a)\, a^{-2}.
$$
For instance, if we set $V=\pa_x$ in the algebra  $\Sym\!(u(2))$, we get the classical result 
$$
\pa_x(x^{-1})=-x^{-2}.
$$

However, the classical Leibniz rule is not valid in the algebra $U(u(2)_h)$. Instead, we use the multiplicativity of the map $\hatt$ (this property can be viewed a 
form of the new Leibniz rule). As was said above, we put $\hatt(a^{-1})=\hatt(a)^{-1}$ by definition, provided the matrix $\hatt(a)$ is invertible in the algebra ${\Mat}(\Ah)$.  
Thus, the map $\hatt$ is well defined on the subset $\tilde{\Ah}\subset \Ah$ consisting of elements $a$ such that the corresponding matrix $\hatt(a)$ is invertible.

Observe that the matrix $\hatt(\rh)$ is invertible. Moreover, it is not difficult to find the matrix $\hatt(\rh^p)=i\hh \Theta(\rh^p)$, $p\in \Bbb{Z}$.
Its entries can be easily computed according to the following formulae (see \cite{GS2})
\be 
\tpa_t(\rh^p)=\frac{-i}{2\hh\rh}((\rh+\hh)^{p+1}+(\rh-\hh)^{p+1}) 
\label{pp} 
\ee
and
\be 
\pa_x(\rh^p)=\frac{x}{\rh}\parh(\rh^p), \qquad \pa_y(\rh^p)=\frac{y}{\rh}\parh(\rh^p),\qquad  \pa_z(\rh^p)=\frac{z}{\rh}\parh(\rh^p), 
\label{ppp} 
\ee
where
$$
\parh(f(\rh))=\frac{f(\rh+\hh)-f(\rh-\hh)}{2\hh}
$$
is the {\em  derivative} in the quantum radius introduced in \cite{GS2}.

For instance, if $p=-1$ we have
$$
\parh(\rh^{-1})=\frac{1}{\hh^2-\rh^2},\qquad \pa_x(\rh^{-1})=\frac{x}{\rh(\hh^2-\rh^2)}\qquad
\tpa_t(\rh^{-1})=\frac{-i}{\hh \rh},\quad\Rightarrow\quad\partial_t(\rh^{-1}) = 0.
$$

Now, we are able to calculate the matrix $\hatt(\rh^p)$, $p \in \Bbb{Z}$. We have
\be
\hatt(\rh^p)=\frac{(\rh+\hh)^{p+1}+(\rh-\hh)^{p+1}}{2\rh} I-\frac{i((\rh+\hh)^{p}-(\rh-\hh)^{p})}{2\rh} A, 
\label{ten} 
\ee
where $A$ is the following matrix
$$
A=
\left(\!\!
\begin{array}{cccc}
 0 & -x\,&-y&-z\\
 x& 0 &-z&y\\
y& z & 0 & -x\\
z&-y&  x& 0
\end{array}
\!\!\right).
$$
This matrix meets the following  relation
$$
A^2=(\hh^2-\rh^2) I-2i\hh\, A.
$$

It can be easily checked that 
$$
\hatt(\rh^p)\hatt(\rh^q)=\hatt(\rh^{p+q}),\qquad \forall\, p,q \in \Bbb{Z}.
$$

If for a given element $a\in \AA$, $a\not=0$ the entries of the matrix $\hatt(a)$ commute with each other, then the inverse matrix $\hatt(a)^{-1}$ with coefficients 
from $\Ah$ can be found in the usual way. So, in this case the matrix $\hatt(a)$ is invertible iff  the determinant of this matrix does not vanish. Note that the determinant 
in this case reads
$$
\det \hatt({a})=-\hh\,(\tpa_t({a})^2+\pa_x({a})^2+\pa_y({a})^2+\pa_z({a})^2)^2.
$$

In general, we can only apply the Gauss method in order to  trigonalize the matrix $\hatt(a)$. It is possible in the skew-field $\Ah$. However, in so doing, we have
to assume a series of elements appearing in computations to be different from 0 in order to invert them.
This condition is not explicit enough. It would be desirable to get a more explicit condition.

Note, that the entries of each matrix $\hatt(t)$, $\hatt(x)$, $\hatt(y)$, $\hatt(z)$ commute with each other. Besides, $\det \hatt(x)=(x^2-\hh^2)^2$ and so on.
So, these matrices are invertible in the algebra $\Ah$. By contrast, the entries of the matrix $\hatt(\rh)$ do not commute with each other but, as we have seen, 
this matrix is invertible. We get the matrix $\hatt(\rh^{-1})$ by putting $p=-1$ in (\ref{ten}).

Now, let us consider a more general example
$$
a=\al_0 t+\al_1 x+\al_2 y+\al_3 z,\quad \al_i\in \K,\quad i=0,1,2,3.
$$
Then the matrix $\hatt(a)$ reads
\be
{\hatt}(a)=
i\hh\left(\!\!
\begin{array}{cccc}\displaystyle
 \al_0 +\frac{a}{i\hh} & -\al_1\,&-\al_2&-\al_3\\
 \al_1& \displaystyle \al_0 +\frac{a}{i\hh}  &-\al_3&\al_2\\
\al_2&\al_3&  \displaystyle \al_0 +\frac{a}{i\hh} & -\al_1\\
\al_3&-\al_2&  \al_1s & \displaystyle \al_0 +\frac{a}{i\hh}
\end{array}
\!\!\right)
\ee

Unfortunately, we have not succeeded  in explicit inverting the matrix $\hatt(\rh-\al_0 t+\al_1 x+\al_2 y+\al_3 z)$ if at least one of the quantities $\al_1$, $\al_2$, $\al_3$
does not vanish. Nevertheless, as follows from our considerations, for all elements $a\in {\cal A}[x^{-1},y^{-1},z^{-1}]$ (noncommutative Laurent polynomials) the matrices
$\hat\Theta (a)$ are invertible.

\section{Quantization with noncommutative configuration space}

Let $\gg$ be a Lie algebra. Then by a quantization of  the commutative algebra $\Sym\!(\gg)$ we mean the enveloping algebra $U(\gg_\h)$ endowed with a map of linear 
spaces $\al_\h: \Sym\!(\gg)\to U(\gg_\h)$ smoothly depending on the quantization parameter $\h$. Such a map enables us to define a new product in the commutative algebra 
$\Sym\!(\gg)$ by means of the following composition of maps:
$$
f\star_h g= \al_h^{-1} ( \al_h(f)\cdot  \al_h(g)), 
$$
where the symbol $\cdot$ stands for the product in the algebra $U(\gg_\h)$. (Consequently, on the algebra $\Sym\!(\gg)$ one can define the corresponding Poisson structure.)

We want such a map $\al_\h$ to be $\gg$-covariant. In the case under consideration  ($\gg=u(2)_\h$) we also want to extend it up to the map ${\A}_0\to \A$ and further up to 
the map ${\Ah}_0\to \Ah$ with keeping the $\gg$-covariance. Hereafter, the notation ${\A}_0$ (resp., ${\Ah}_0$) stands for the classical limit  of the algebra $\A$ (resp., $\Ah$) 
as $\h\to 0$. A way of constructing such a map ${\A}_0\to \A$ is presented in the paper \cite{GS2}. It is based on the Kostant's theorem saying that for a simple Lie algebra $\gg$ the enveloping algebra $U(\gg)$ is free over its center.

In the framework of our approach with any function $f(t, r)\in \K(t,r)$ we associate the element $f(t, \rh)$. On the next step with a fraction $f/g$, $f,g \in \AA_0$ we associate  
the element\footnote{Here we associate the right fraction with the element $f/g$. However, in principle, we can use left fractions or their combinations.} $\al_h(f)(\al_\h(g))^{-1}$. Consequently, we have a map $\al_h: \Ah_0\to \Ah$ (we keep for the extended map the same notation).

However, we are aiming at the quantization of  the algebra of differential operators and the de Rham complex, related to the algebra $\Ah$. We call the corresponding 
procedure {\em  the quantization with noncommutative configuration space}. Emphasize that in the canonical quantization, giving rise to operators acting on a Hilbert space,
the noncommutativity appears only at the level of the phase space.

In order to quantize a differential operator $P$ with coefficients from $\Ah_0$ we first quantize these coefficients by means of the method described in \cite{GS2}. Then
we convert the partial derivatives in the operator $P$ into the quantum partial derivatives. Thus, we get an operator acting on (a subset of) the algebra $\Ah$. We denote 
it $\al_h(P)$. In general, $\al_h(P(f))\not=\al_h(P)(\al_\h(f))$, $f\in \Ah_0$.

Now,we describe the corresponding quantization of the de Rham complex. However, first we construct an analog of de Rham complex on the algebra $\Ah$ (more precisely, 
on the aforementioned subset).

Let us introduce new generators $dt$, $dx$, $dy$, $dz$ called {\em  pure differentials}, which are assumed to anti-commute with each other. Let $\wedge^k$, $k=1,2,3,4$  
be the homogeneous components of the skew-symmetric algebra generated by the pure differentials. Then our de Rham operator is defined by
$$
d:\wedge^k\ot \Ah\to \wedge^{k+1}\ot \Ah,\quad d(\omega\ot f)=\sum_{u\in \{t,x,y,z\}} \omega\wedge du\ot \pa_u(f),\quad \forall\, \omega\in \wedge^k, \, f\in \Ah.
$$

We leave checking the property $d^2=0$ to the reader. Similarly to the classical case, we assume that the first and the second quantum partial derivatives of the element $f$
in the formula above are well defined.

Now, we are able to extend the map $\al_h$ to the terms of the de Rham complex over the algebra $\Ah_0$. First, to each element from the space $\wedge^k$ generated
by $dt$, $dx$, $dy$, $dz$, where $t, x, y, z\in \Sym\!(u(2))$, we assign the same element but with $t$, $x$, $y$, $z$ treated as elements of $U(u(2)_\h)$. Then, with any 
element
$$
\omega\ot f,\quad \omega\in \wedge^k, \quad f\in \Ah_0
$$
we associate the element $\al_h(\omega) \ot \al_h(f)$ where $\al_h(\omega)$ is as above. Note that, in general, $\al_h(d \omega)\not=d \al_h(\omega)$.

\begin{remark}
\rm There is known another Lie algebra $\gg$ such that its enveloping algebra permits defining quantum derivatives. This Lie algebra is defined by the
following table of Lie brackets:
\be 
[x,y]=[y,z]=[z,x]=0,\quad [x, t]= \kappa^{-1} x,\quad [y, t]=\kappa^{-1} y,\quad
[z, t]=\kappa^{-1} z. 
\label{Lie} 
\ee
Its enveloping algebra is called $\kappa$-Minkowski space \cite{MR}.

The ``quantum partial derivatives" on this algebra are defined on the generators in the usual way, but
the Leibniz rule is modified according to the coproduct
$$
\De(\pa_t)=\pa_t\ot 1+1 \ot \pa_t,\quad \De(\pa_u)=\pa_u\ot 1+ \exp(-\kappa^{-1}\pa_t)\ot \pa_u,\quad u\in\{x,y,z\}.
$$

It is possible to quantize the algebra $\Sym\!(\gg)$ and the corresponding algebra of differential operators in the manner discussed at the beginning of this section.
However, the corresponding quantum calculus is much less rich. In, particular, this algebra does not enable one to quantize the radius.

One more algebra, which can be considered as a NC configuration space, stems from Moyal quantization of a constant Poisson bracket. For such an algebra
the corresponding ``phase space" can be constructed with the help of the usual partial derivatives. Nevertheless, there are no such structures  covariant with respect to
the group $SO(3)$.
\end{remark}

\section{Noncommutative Dirac monopole}

Let us quantize the Maxwell system on the Minkowski  space by the method presented in the previous sections.

As is known, the Maxwell system consists of 4 equations. The first couple of these equations is (we put $c=1$)
$$
\di \H=0,\quad {\rm rot}\E+\pa_t \H=0,
$$
where  ${\E}=(E_1, E_2, E_3)$ and ${\H}=(H_1, H_2, H_3)$ are the vectors of electric and magnetic fields respectively. 

The second pair of the Maxwell equations in vacuum is
$$
\di \E=0,\quad {\rm rot}\H- \pa_t\E=0.
$$

Let us consider a particular case of this system, giving rise to the Dirac monopole, i.e.  we assume $\E$ to be zero, and consequently
$\H$ is assumed to be stationary. Then, we get the following system for the magnetic field
\be
\di \H=0,\quad {\rm rot}\H=(0,0,0). 
\label{mon} 
\ee

We look for a spherical symmetric solution of the system (\ref{mon}), that is, ${\H}=f(r)(x,y,z)$, where $f(r)$ is a function in the radius $r$ to be defined. From the first 
equation of this system we get the equation on $f$:   
$$ 
3f+r\,\frac{df}{dr}=0.
$$ 
This equation has the general solution $f(r)=g\, r^{-3}$  where $g$ is a constant which is assumed to be real. However, the field $\H=\frac{g}{r^3}(x,y,z)$ is a solution of 
the equation $\di \H=0$ only on the set $\R^3\setminus (0,0,0)$, whereas, on the whole space $\R^3$ this field meets the equation
\be 
\di {\H}= 4\,g \pi \de(r),  
\label{div}   
\ee
where $\de(r)$ is the delta-function on the space $\R^3$ located at the point $(0,0,0)$.

As for the second equation of  (\ref{mon}), it is satisfied by ${\H}=f(r)(x,y,z)$  with any rational function $f(r)$.

Now, let us  apply our quantization to this model. Since the Maxwell system and consequently, the system (\ref{mon}) consists of operators with constant coefficients, its 
quantization in our sense is somewhat easy. We only have to replace all partial derivatives by their quantum counterparts. Again, let us try to solve the system (\ref{mon}) 
by looking for a spherical symmetric solution ${\H}=(x,y,z)f(\rh)$, where $\di$ and ${\rm rot}$ are expressed in terms of the quantum partial derivatives via the classical formulae.

First, consider the second equation. Let us compute the first coordinate of the vector ${\rm rot} \H$. We have
\begin{eqnarray*}
({\rm rot}\H)_x &=&\pa_y(zf(\rh))-\pa_z(yf(\rh))=(z\pa_y+\frac{h}{2}\pa_x-y\pa_z+\frac{h}{2}\pa_x)f(\rh)\\
&=& (z\frac{y}{\rh}-y\frac{z}{\rh}+h\frac{x}{\rh})\parh (f(\rh))=0
\end{eqnarray*} 
for any rational function $f$. Other components of the vector ${\rm rot}\H$ vanish as well. Thus, similarly to the classical case the equation ${\rm rot}\H=(0,0,0)$ is fulfilled 
with any rational $f$.

Now, consider the  equation $\di\, \H=0$. We have
\begin{eqnarray*}
&&\pa_x(xf(\rh))+\pa_y(yf(\rh))+\pa_z(zf(\rh))=(x\pa_x+y\pa_y+z\pa_z+\frac{3h}{2}\tpa_t)(f(\rh))\\
&&=(\frac{x^2+y^2+z^2}{\rh}\parh+3i\hh \tpa_t)(f(\rh))=\frac{\rh^2-\hh^2}{\rh}\frac{f(\rh+\hh)-f(\rh-\hh)}{2\hh}\\
&&+3\frac{f(\rh+\hh)(\rh+\hh)+ f(\rh-\hh)(\rh-\hh)}{2\rh}=0.
\end{eqnarray*}
Here, we used the formula
$$
\tpa_t(f(\rh))=\frac{f(\rh+\hh)(\rh+\hh)+f(\rh-\hh)(\rh-\hh)}{2i\hh \rh}
$$
which results from (\ref{pp}).

This entails the following relation on the function $f$:
$$
(\rh^2+2\hh^2+3\rh\hh)f(\rh+\hh)=(\rh^2+2\hh^2-3\rh\hh)f(\rh-\hh).
$$
Having introduced a new function $\psi(\rh)=\rh f(\rh)$ we get the following difference equation:
$$
(\rh+2\hh)\psi(\rh+\hh)=(\rh-2\hh)\psi(\rh-\hh).
$$
It can be checked that the function
$$
\psi(\rh)= \frac{g}{\rh^2-\hh^2},\quad g\in \K
$$
is a solution of this equation ($g$ is constant for a fixed $\hh$ but in principle it can depend on $\hh$). 

Consequently, the NC field
\be 
\H=\frac{g}{\rh (\rh^2-\hh^2)} (x,y,z) 
\label{mag} 
\ee
meets the system (\ref{mon}). Note that its components are elements of the algebra $\AA$ (here we do not need the corresponding skew-field).

We call this solution {\em the NC Dirac monopole}. Emphasize that for $\hh\to 0$ we retrieve the classical Dirac monopole.

\begin{remark} 
\rm 
Another solution can be found by using the Fourier transform. Omitting details of calculation, we only present the final result:
$$
\psi(\rh)=C(\hh)(\de(\rh+\hh)-\de(\rh-\hh)).
$$
Here $C(\hh)$ is an arbitrary function in $\hh$. The corresponding function $f$ reads
$$
f(\rh)=B(\hh)(\de(\rh+\hh)+\de(\rh-\hh)),
$$
where $B(\hh)=-C(\hh)\hh^{-1}$.

However, the $\de$-function entering this formula does not belong to the class $\K(t,\rh)$ and we disregard such a solution.
\end{remark}

Nevertheless,  similarly to the classical case, it is reasonable to assume that that the equation $\di \H=0$ should be corrected on the whole algebra $\A$.
Using a heuristical reasoning, we can argue that for the NC Dirac monopole (\ref{mag}) the relation (\ref{div}) holds. We follow the guidelines of \cite{GS2} 
(see section 6).

In the classical case the relation (\ref{div}) considered in the framework of distribution theory means that
\be 
\int \left({x}{f(r)}\pa_x(\varphi)+{y}{f(r)}\pa_y(\varphi)+{x}{f(r)}\pa_z(\varphi)\right)dx\, dy\, dz=-4\pi \varphi(0), 
\label{int} 
\ee
where  the integral is treated in the sense of the principle value. Also, a test function $\varphi$ is assumed to belong to $\AA_0$ and to tend to 0 as $r\to \infty$.
Since the integral in (\ref{int}) is $SO(3)$-invariant, it is possible to suppose that the function $\varphi$ depends on $r$ only.

Passing to the the algebra $\AA$ we have to define a quantum counterpart of the integral (\ref{int}). In (\ref{int}) we replace $f(r)$ for $f(\rh)=(\rh(\rh^2-\hh^2))^{-1}$
and  define the quantum analog of the integral (\ref{int}) as an $SO(3)$-invariant map
\be 
\AA\to \K,\qquad \varphi\mapsto \int \varphi\,dx\,dy\, dz. 
\label{intt} 
\ee
Consequently, the map (\ref{intt}) kills all components of the decomposition $\AA=\oplus (\K(t,\rh)\ot V_i)$, where $V_i$ run over all  irreducible $su(2)$-modules of 
integer spin, except for the trivial one $\K(t,\rh) $. Thus, the relation (\ref{int}) makes sense if we assume that  $\varphi(\rh)\in \K(t,\rh)$ and $\varphi$ is regular on $\R$. 
As for the normalization of the integral on the trivial component $\K(t,\rh)$, its follows from the first formula (\ref{bel}) below (in fact, it is a definition).

Then the left hand side of (\ref{int}) can be presented as follows
$$
\lim_{\epsilon\to 0\atop R\to\infty} \int_{B(R)\setminus B(\epsilon)} \left(\frac{x}{\rh(\rh^2-\hh^2)} \frac{x}{\rh} (\parh \varphi)+
\frac{y}{\rh(\rh^2-\hh^2)} \frac{y}{\rh} (\parh \varphi)+\frac{z}{\rh(\rh^2-\hh^2)} \frac{z}{\rh} (\parh \varphi)\right)dx dy dz
$$ 
$$
 =\lim_{\epsilon\to 0\atop R\to\infty} \int_{B(R)\setminus B(\epsilon)}\frac{x^2+y^2+z^2}{\rh^2(\rh^2-\hh^2)} (\parh \varphi)\, \rh^2\, d\, \rh\, 
d\,\Omega=4\pi \lim_{\epsilon\to 0\atop R\to\infty} \int_{B(R)\setminus B(\epsilon)}(\parh \varphi)d\,\rh
$$
$$
=4\pi\lim_{\epsilon\to 0\atop R\to\infty}(\varphi(R)-\varphi(\epsilon))=-4\pi\,\varphi(0).
$$
Here $B(a)$ stands for a ball with the center at 0 of radius $a\geq 0$, and similarly to \cite{GS2}, we assume that analogs of the classical formulae
\be 
dx\, dy\, dz=\rh^2\, d\,\rh\, d\,\Omega,\qquad \int_a^b \parh(\varphi(\rh))\, d\,\rh=\varphi(b)-\varphi(a), 
\label{bel} 
\ee
are valid, where $d\,\Omega$ is the volume form of the unit sphere. Undeer this assumption, the NC Dirac magnetic monopole (\ref{mag}) meets the 
relation (\ref{div}).

Note that in the classical case there exist vector potentials for the magnetic field $\H$, i.e. a vector-functions $\AAA=(A_1,A_2,A_3)$ such that ${\rm rot}\AAA=\H$. 
These vector potentials are parameterized by a unit vector $\nn=(n_1,n_2,n_3)$ and have the form
$$
\AAA=\frac{g}{r}\frac{(x,y,z)\times \nn}{(r-(x,y,z)\cdot \nn)},
$$
where $\times$ is the vector product of two vectors, and $\cdot$ is their scalar product. Each of these vector-potentials
is singular on a half-line.

The image $\al_h(\AAA)=(\al_h(A_1),\al_h(A_2),\al_h(A_3))$ of any such a vector-potential its image is well-defined: its components belong to the algebra $\Ah$.
However, we are not able to compute the vector ${\rm rot} {\al_\h(\AAA)}$ explicitly, since we do not know how to invert the matrix $\hatt(\rh-(n_1 x+n_2 y+n_3 z))$ 
(see section 3).

Also, note that in the classical setting the Maxwell system can be realized in terms of differential forms. By using the quantizing  map $\al_h$ extended to the space
of differential forms over  the algebra $\Ah_0$ we can quantize the Maxwell system in this way. However, such a quantization leads to the same system on the ``NC fields" 
$\E$ and $\H$ and, consequently, to the same NC Dirac monopole. We leave details to the reader.

\begin{remark}\rm 
In the classical setting the second pair of Maxwell equations (and consquently, the second equation of (\ref{mon})) arises from a metric. We do not consider any metric in the
present paper and impose the equations by transferring them from the classical algebra. By chance, the equation ${\rm rot} \H$ is fulfilled for any $f$. We plan to consider 
models with metric in our subsequent publications.
\end{remark}


\begin{thebibliography}{999}

\bibitem{GPS1} Gurevich D., Pyatov P., Saponov P.:
Cayley-Hamilton theorem for quantum matrix algebras of $GL(m|n)$ type.  St. Petersburg Math. J. 17(1), 119–135 (2006).

\bibitem{GPS2} Gurevich D., Pyatov P., Saponov P.:
Braided Weyl algebras and differential calculus on $U(u(2))$.  J. Geom. Phys. 62(5), 1175–1188 (2012).

\bibitem{GS1} Gurevich D., Saponov P.: Braided algebras and their applications to Noncommutative Geometry.  Adv. in Appl. Math. 51(2), 228–253 (2013).

\bibitem{GS2} Gurevich D., Saponov P.: Noncommutative Geometry and dynamical models on $U(u(2))$ background. J. Generalized Lie Theory Appl 9(1) (2015).

\bibitem{GS3} Gurevich D., Saponov P.: Derivatives in noncommutative calculus and deformation property of quantum algebras.  J. Noncommuttative Geometry, to be published

\bibitem{MR} Majid Sh., Ruegg H.:  Bicrossproduct structure of $\kappa$-Poincaré group and non-commutative geometry.
Phys. Lett. B 334(3-4), 348–354 (1994).

\end{thebibliography}
\end{document}